\documentclass[11pt]{article}

\usepackage[centertags]{amsmath}
\usepackage[normal]{subfigure}
\usepackage{wrapfig}
\usepackage{verbatim}
\usepackage{amsfonts}
\usepackage{amssymb}
\usepackage{amsthm}
\usepackage{amsmath}
\usepackage{graphicx}
\usepackage{epsfig}
\usepackage{amsmath}
\usepackage{color}
\usepackage{multicol}
\usepackage{amssymb,amsmath}
\usepackage{latexsym}
\usepackage[normalem]{ulem}
\usepackage{cancel}

\usepackage{bm}

\usepackage{amsmath, amsthm, amssymb}

\newtheorem{theorem}{Theorem}[section]

\newtheorem{e-proposition}[theorem]{Proposition}

\newtheorem{e-definition}[theorem]{Definition\rm}

\DeclareTextFontCommand{\rm}{\rmfamily}
 \usepackage{color}
\definecolor{dred}{rgb}{0.92,0,0}
\definecolor{dgreen}{rgb}{0,0.92,0}
\definecolor{dblue}{rgb}{0,0,0.92}
\definecolor{dyellow}{rgb}{0.95,0.95,0}

\newenvironment{prooff}{{\it Proof :}}{\hfill\rule{2mm}{2mm}\vskip3mm
\par}

\newcommand{\ds}{\displaystyle}
\newcommand{\pa}{\partial}

\def\Ecal{{\bf \mathcal{E} }}

\def\nuv{\mbox{\boldmath$\nu$}}
\renewcommand{\div}{{\rm div}\,}

\newcommand{\curlv}{{\bf curl}\,}

\newcommand{\gradv}{{\bf grad}\,}

\def\vvec{{\bf v}}
\def\xvec{{\bf x}}

\newcommand{\Evec}{\mathbf{E}}
\newcommand{\Bvec}{\mathbf{B}}
\newcommand{\Fvec}{\mathbf{F}}
\newcommand{\Jvec}{\mathbf{J}}

\newcommand{\partialderiv}[2]{\frac{\partial #1}{\partial #2}}

\newcommand{\bea}{\begin{eqnarray}}
\newcommand{\eea}{\end{eqnarray}}
\newcommand{\fracp}[2]{\frac{\partial #1}{\partial #2}}
\def\R{\mathbb{R}}
\newcommand{\p}[1]{\left(#1\right)}

\def\N{\mathbb{N}}
\newcommand{\khat}{\hat{\bold{e}}_z}
\newcommand{\curlvp}{{\bf curl_{\perp}}}
\newcommand{\curlsp}{{\rm curl_{\perp}}}
\newcommand{\pp}[1]{\left[#1\right]}

\newcommand*\diff{\mathop{}\!\mathrm{d}}
\newcommand{\divp}{{\rm div_{\perp}}} 
\newcommand{\be}{\begin{equation}}
\newcommand{\ee}{\end{equation}}
\def\nuv{\mbox{\boldmath$\nu$}}
\def\tauv{\mbox{\boldmath$\tau$}}

\newtheorem{Lem}{Lemma}[subsection]

\title{A hierarchy of reduced models to approximate Vlasov-Maxwell equations for slow time variations}
\author{
F. Assous\footnote{Ariel University, 40700 Ariel, Israel.} ,
Y.Furman \footnote{Ariel University, 40700 Ariel, Israel.}
}
\date{}

\begin{document}
\maketitle
\begin{abstract}
We introduce a new family of paraxial asymptotic models that approximate the Vlasov-Maxwell equations in non-relativistic cases. 
This formulation is $n$-th order accurate in a parameter $\eta$, which denotes the ratio between the characteristic velocity of the beam and the speed of light. This family of models is interesting, first because it is simpler than the complete Vlasov-Maxwell equation, then because it allows us to choose the model complexity according to the expected accuracy.
\end{abstract}

{\bf\noindent keywords}: Vlasov-Maxwell equations;  asymptotic analysis;  paraxial model.
\section{Introduction}
\label{intro}

\noindent Charged particle beams are very useful in a variety of scientific and technological applications. After the discovery that both magnetic and electric fields can act as lenses for electron rays, the field experienced rapid development, with industrial applications such as welding \cite{weld}, micromachining and lithography \cite{lith}, thermonuclear fusion \cite{fusion}, etc. More recent developments use intense electron beams as electromagnetic radiation sources, like the gyrotron or the free-electron laser (see for instance \cite{gyrotron},\cite{tran}). More details can be found in \cite{Laws88} and \cite{reiser2008book}. Hence, there is a great interest in mathematical and numerical modeling of these phenomena.\\

\noindent Considering non-collisional beams, a well-accepted method for describing the transport of bunches of particles is the Vlasov equation (\cite{vlasov}, \cite{BiLa85}). Since the particles are electrically charged, the force field which governs their movement is the Lorentz force, which in turn depends on both the electric and magnetic fields, which are solutions to the well-known Maxwell equations \cite{AsCL18}. This set of equations coupled together is known as the time-dependent Vlasov-Maxwell system of equations.\\

\noindent However, the numerical solution of this model, which is unavoidable in many situations, \cite{ADHRS93}, \cite{AsDS92}, requires a large computational effort, usually based on a combination of finite elements or finite volume discretisation with particle-in-cell methods. Therefore, whenever possible, it is worthwhile to take into account the particular details of the problem in order to derive approximate models, leading to cheaper simulations (see \cite{DeRa93}, \cite{LaMR94}, \cite{MoMW96}, \cite{RaSo96}).\\

\noindent Following the principle exposed in \cite{LaMR94}, our approach relies on the introduction of a moving frame, which travels along the optical axis at a given velocity. Many noticeable research works have been done in this field: in the case of high-energy, ultra-relativistic short beams, Laval \textit{et al.} \cite{LaMR94}  derived a paraxial approximation of the  Vlasov-Maxwell equations, by introducing a moving frame, which travels along the optical axis at the speed of light $c$.\\

\noindent  This idea of changing variables to follow the moving frame is not new, and can be found elsewhere, for instance in \cite{boyd}, \cite{slinker}. Similar work has been done for a laminar beam case in \cite{nouri}.  A different paraxial model was also derived for the case of high-energy short beams \cite{AsCh12}, and was typically related to free-electron lasers or particle accelerators. This work takes into account the specific geometrical features of the devices, leading thus to a somewhat different dimensional analysis. Numerical applications were also proposed in \cite{AsTs09}, whereas comparison methods of these  models, based on data mining techniques, have been proposed in \cite{AsCh11}.\\

\noindent The aim of this paper is to derive a a new family of paraxial asymptotic models that approximate the Vlasov-Maxwell equations in non-relativistic cases. Section \ref{V-Mmodel} gives a short overview of the equations and the change of variables to the beam frame. The scaling of the equations is presented Section \ref{Scaling}, whereas we propose, in Section \ref{Asymptotic}, the asymptotic expansion of the relevant parameters to derive a new family of paraxial model. Finally, the resulting paraxial models, that allow us to choose the model complexity according to the expected accuracy, are given in Section \ref{ParaxialModel}.

\section{The Vlasov-Maxwell model}
\label{V-Mmodel}
Consider a beam of charged particles with a mass $m$ and a charge $q$ moving in a perfectly
conducting cylindrical tube, whose axis is constituted by the $z$-axis. We denote by  $\Omega$ the transverse section of boundary $\Gamma$,  $\nu=(\nu_x,\nu_y,0)$ denoting the unit exterior normal to the tube. We suppose that an external magnetic field $\Bvec^e$ confines the beam in a neighbourhood of the $z$-axis which may be therefore chosen as the optical axis of the beam. Let $\xvec=(x,y,z)$ be the position of the particle and $\vvec=(v_{x}, v_{y}, v_{z})$ its velocity. We assume that the beam is non relativistic and non collisional so that its distribution function $f=f(\xvec,\vvec,t)$ in the phase space $(\xvec,\vvec)$ is a solution to the Vlasov
equation
\begin{equation}
\label{eq:VlasovOriginal_Eq}
   \frac{\partial f}{\partial t}+\textbf{v}\cdot\gradv_\xvec f +\frac{1}{m}\Fvec\cdot\gradv_\vvec f=0\,.
\end{equation}
Above, $\Fvec=q(\Evec+\vvec\times \Bvec)$ denotes the electromagnetic force acting on the particles. The electric field 
$\Evec=\Evec(\xvec,t)$ and the magnetic field $\Bvec=\Bvec(\xvec,t)$ are solutions to Maxwell's equations
\bea
\frac{1}{c^2}\fracp{\Evec}{t}-\curlv\Bvec&=&-\mu_0\Jvec\\
\fracp{\Bvec}{t}+\curlv\Evec&=&\mathbf{0}\\
\div\Evec&=&\frac{\rho}{\varepsilon_0}\\
\div\Bvec&=&0
\eea
where the charge and the current density $\rho(\xvec,t)$ and $\Jvec(\xvec,t)$ are obtained from the distribution function $f(\xvec,\vvec,t)$ with
\begin{equation}
\label{eq:rho_Eq}
 \rho(\xvec,t) = q \int_{\R^3_\vvec} f d\vvec,\qquad \qquad \Jvec (\xvec,t)= q \int_{\R^3_\vvec} \vvec f d\vvec.
\end{equation}
Now, we introduce a parameter $0 < \beta < 1$, and we consider that the particle longitudinal velocity $v_z$ satisfies $v_z\simeq \beta c$ for any particle in the beam. Hence, we rewrite the Vlasov-Maxwell equations in a frame which moves along $z$-axis with the velocity $\beta c$, i.e. a fraction of the light velocity. For this purpose,  we set $\zeta=\beta c t- z$, $v_\zeta=\beta c-v_z$ and we perform the change of variables 
$(x,y,z,v_x,v_y,v_z,t) \rightarrow (x,y,\zeta,v_x,v_y,v_{\zeta},t)$, so that
\begin{equation}
\p{\fracp{}{z},\fracp{}{v_z},\fracp{}{t}}\to\p{-\fracp{}{\zeta},-\fracp{}{v_\zeta},\fracp{}{t}+\beta c\fracp{}{\zeta}}\,.
\end{equation}
It is also convenient to introduce the transverse quantities\\

\begin{center}
 $\textbf{x}_\perp=(\emph{x},\emph{y}) $,
$\textbf{v}_\perp=(\emph{$v_x$},\emph{$v_y$})$
\end{center}
and to define the transverse operators \\
\begin{center}
 $\textbf{grad}_\perp \varphi=(\ds\frac{\partial \varphi}{\partial x},\frac{\partial \varphi}{\partial
 y})$,
\quad
$
 \textbf{curl}_\perp \varphi=(\ds\frac{\partial \varphi}{\partial y},-\frac{\partial \varphi}{\partial x})$,
 \quad
 $
 \Delta_\bot  \varphi=\ds\frac{\partial^2 \varphi}{\partial x^2}+\frac{\partial^2 \varphi}{\partial y^2}\,,
 $ \\[0.2cm]
\end{center}
where $\varphi=\varphi(x,y)$ is a scalar function. Similarly, for $\textbf{A}_\perp=(A_x,A_y)$ denoting a transverse vector field, we set \\
\begin{center}
   $\div_{\perp}{\textbf{A}_\perp}$=$\ds\partialderiv{A_x}{x}$+$\ds\partialderiv{A_y}{y}$,
   \quad$\curlsp{\textbf{A}_\perp}$=$\ds\partialderiv{A_y}{x}$-$\ds\partialderiv{A_x}{y}$. \\[0.2cm]
\end{center}
We define $\textbf{A}_{\bot}\times \textbf{e}_z=(A_y,-A_x)$ and we readily get the following identities
\begin{equation}
\label{relbas}
 \div_\bot(\textbf{A}_\bot\times\textbf{e}_z)=\curlsp{A},\,\,\quad
  \curlsp(\textbf{A}_\bot\times\textbf{e}_z)=-\div_{\bot}{A},\,\,\quad
  \curlsp\textbf{curl}_\perp \varphi=-\Delta_\perp \varphi.
\end{equation}

\noindent Moreover, denoting by $\tauv=(-\nu_y,\nu_x)$ the unit tangent along $\Gamma$, we have the relation
$\textbf{curl}_\perp \varphi\cdot \tauv= -\ds\frac{\pa \varphi}{\pa \nuv}$.\\

\noindent  Using the above notations, the Vlasov equation in the new variables can be written as
\begin{equation}
\label{VlasovNV}
\fracp{f}{t}+\vvec_\perp\cdot\textbf{grad}_\perp f+v_\zeta\fracp{f}{\zeta}+\frac{1}{m}\Fvec_\perp\cdot\textbf{grad}_{\vvec_\perp}f-\frac{F_z}{m}\fracp{f}{v_\zeta}=0\,.
\end{equation}

\noindent  Additionally, setting $\Ecal_\perp=\p{E_x-\beta cB_y,E_y+\beta cB_x}$ and $J_\zeta=\rho\beta c-J_z=q\int_{\mathbb{R}^{3}_\vvec} v_\zeta f \diff \vvec$, we obtain the following expressions for Maxwell's equations. First, Gauss's law that takes the form:
\begin{equation}
\label{GaussE1}
\divp\Evec_\perp-\fracp{E_z}{\zeta}=\frac{\rho}{\varepsilon_0}\,,
\end{equation}
and Gauss's law for magnetism is expressed
\begin{equation}
\divp\Bvec_\perp-\fracp{B_z}{\zeta}=0\label{eq:eq4}\,.
\end{equation}
In the same way, Ampere's law can be written 
\begin{eqnarray}
\perp&:&\frac{1}{c^2}\fracp{\Evec_\perp}{t}+\frac{1}{\beta c}\fracp{}{\zeta}\p{\Ecal_\perp-\p{1-\beta^2}\Evec_\perp}-\curlvp B_z=-\mu_0\Jvec_\perp \label{Amp11}\,,\\
\zeta&:&\frac{1}{c^2}\fracp{E_z}{t}+\frac{1}{\beta c}\divp\p{\Ecal_\perp-\p{1-\beta^2}\Evec_\perp}=\mu_0 J_\zeta \label{Amp12}\,,
\end{eqnarray}
and Faraday's law becomes
\begin{eqnarray}
\perp&:&\fracp{\Bvec_\perp}{t}+\fracp{}{\zeta}\p{\Ecal_\perp\times\khat}+\curlvp E_z=\mathbf{0}\label{eq:eq2}\,,\\
\zeta&:&\fracp{B_z}{t}+\curlsp\Ecal_\perp=0\label{eq:eq3}\,.
\end{eqnarray}
Finally, the electromagnetic force becomes
\begin{eqnarray}
\label{FocesNV}
\Fvec_\perp&=&q\p{\Ecal_\perp+\p{B_z\vvec_\perp+v_\zeta\Bvec_\perp}\times\khat}\,,\\
F_z&=&q\p{E_z+\vvec_\perp\cdot\p{\Bvec_\perp\times\khat}}\,.
\end{eqnarray}
Let us formulate now the boundary conditions. Assuming that the particles remain inside a fixed domain $\Omega\times\p{0,Z}$ in the beam frame,  it means that 
$f=0$ on the boundary. For the initial conditions, we simply assume that the initial distribution of particles is a known function which satisfies the boundary conditions $f_{\mid {t=0}}=f_0$.\\
Regarding the electromagnetic fields, the surface of the tube being a perfect conductor, the tangential components of the electric field vanish, for ${\xvec_\perp\in\Gamma,~\zeta\in\p{0,Z}}$ and we have 
$$
\Evec_\perp\cdot\tauv=0, \qquad E_z=0\,.
$$
For the artificial boundary $\zeta=0$, assuming there is no external electric field, and that the static electromagnetic fields that exist ahead of the beam cannot be modified by the electromagnetic waves generated by the beam, we have for 
${\xvec_\perp\in{\Omega},~\zeta=0}$:
$$
\Evec=0, \qquad \Bvec=\Bvec^e, \mbox{ where } \Bvec^e \mbox{ denotes a given external field}\,.
$$
We also assume given initial conditions $\Evec|_{t=0}=\Evec_0, \Bvec|_{t=0}=\Bvec_0$, where $\Evec_0$ and $\Bvec_0$ satisfy both Maxwell's equations and the boundary conditions specified above.\\

\noindent  Let us note some important consequences, for the sequel, of these boundary conditions. Taking the inner product of $\Ecal_\perp$ and $\tauv$ for ${\xvec_\perp\in\Gamma,~\zeta\in\p{0,Z}}$ we get:
\be
\label{BC-Ecal}
\Ecal_\perp\cdot\tauv=\beta c\Bvec_\perp\cdot\nuv\,.
\ee
Next, taking the dot product of (\ref{eq:eq2}) by $\nuv$ and using the definition of $\curlvp$, one obtains for $\xvec_\perp\in\Gamma,~\zeta\in\p{0,Z}$:
\be
\label{eq:eq05}
\p{\fracp{}{t}+\beta c\fracp{}{\zeta}}\p{\Bvec_\perp\cdot\nuv}=0\,.
\ee
Similarly, integrating (\ref{eq:eq3}) over $\Omega$ and applying Green's theorem for $\zeta\in\p{0,Z}$ we get:
\be
\label{eq:eq051}
\int_\Omega\fracp{B_z}{t}\diff\xvec_\perp+\beta c\oint\limits_\Gamma\Bvec_\perp\cdot\nuv\diff l=0\,.
\ee
In the same spirit as above, we obtain using (\ref{eq:eq4}), for $\zeta\in (0,Z)$:
\be
\label{eq:eq052}
\int_\Omega\p{\fracp{}{t}+\beta c\fracp{}{\zeta}}B_z\diff\xvec_\perp=0\,.
\ee

\section{A scaling of the equations}
\label{Scaling}

The second step to derive the paraxial model is to introduce an {\em ad hoc} scaling of the equations. Assuming that we deal with a short beam, we introduce a scaling of the equations by handling the following properties of the beam:
\begin{enumerate}
\item  The beam dimension is small compared to the longitudinal length $L$ of the device;
\item The transverse particle velocities $\textbf{v}_{\bot}$ are comparable to $v_\zeta$, so we have $v_\zeta\simeq\vvec_\perp \ll v_z\simeq\beta c$.
\end{enumerate}
\noindent  Thus, we introduce the two characteristic quantities:
\begin{enumerate}
\item $l$, the characteristic dimension of the beam,
\item  $\overline{v}$,  the characteristic velocity of the particles.
\end{enumerate}

\noindent Note that, in contrast to the case described in \cite{LaMR94}, \cite{AsTs09} or \cite{AsCh14}, we did not require here the longitudinal particle velocities $v_z$ to be necessary close to the light velocity $c$, since we consider a non-relativistic case. For this reason, we set $v_z \simeq \beta c, 0 < \beta < 1$, which allows us to play on the value of the parameter $\beta$.\\

\noindent Now,  defining a small parameter $\eta$ and a characteristic time $T$ with 
\begin{equation}\label{eta}
\eta \equiv \frac{\overline{v}}{c}\ll 1, \qquad T=\frac{l}{\bar{v}}\,,
\end{equation}
we can write:
\be
\begin{array}{l}
x=lx',~y=ly',~\zeta=l\zeta',~t=Tt',~v_x=\bar{v}v_x',~v_y=\bar{v}v_y',~v_\zeta=\bar{v}v_\zeta'
\end{array}
\ee
where the primes represent dimensionless quantities. Using the physical units of the physical quantities and based on the Vlasov-Maxwell equations, one can introduce the following scaling factors: For the electric field one can define $\bar{E}=\ds\frac{m\bar{v}^2}{ql}$, so that from Gauss's law, one can set $\bar{\rho}=\ds\frac{\varepsilon_0 m\bar{v}^2}{ql^2}$. 
From the definition of $\rho$ we get $\bar{f}=\ds\frac{\varepsilon_0m}{q^2l^2\bar{v}}$. Similarly, using the physical units of the other quantities we obtain that $\bar{J}=\ds\frac{\varepsilon_0mc\bar{v}^2}{ql^2}, \bar{F}=\ds\frac{m\bar{v}^2}{l}$ and $\bar{B}=\ds\frac{m\bar{v}^2}{qcl}$.
This allows us to write $f(\xvec_\perp,\zeta,\vvec_\perp,v_\zeta,t)=\bar{f}f'(\xvec_\perp',\zeta',\vvec_\perp',v_\zeta',t')$, 
$\Evec(\xvec_\perp,\zeta,t)=\bar{E}\Evec'(\xvec_\perp',\zeta',t')$, 
$\Bvec(\xvec_\perp,\zeta,t)=\bar{B}\Bvec'(\xvec_\perp',\zeta',t')$ and $\Fvec(\xvec_\perp,\zeta,\vvec_\perp,v_\zeta,t)=\bar{F}\Fvec'(\xvec_\perp',\zeta',\vvec_\perp',v_\zeta',t')$.\\

\noindent Now, defining $\rho'=\ds\int_{\mathbb{R}^{3}_\vvec}f'\diff \vvec'$ and $\Jvec'=\ds\int_{\mathbb{R}^{3}_\vvec}\vvec'f'\diff \vvec'$, it is convenient to introduce $\rho=\bar{\rho}\rho'$ for the charge density, and
$\Jvec_\perp=\bar{J}\eta\Jvec_\perp', J_\zeta=\bar{J}\eta J_\zeta'$ for the current density.\\

\noindent Hence, we are able to write down the Vlasov-Maxwell equations using these dimensionless variables. Dropping the primes for simplicity, the Vlasov equation in dimensionless variables is simply
\begin{eqnarray}
\fracp{f}{t}+\vvec_\perp\cdot\textbf{grad}_\perp f+v_\zeta\fracp{f}{\zeta}+\Fvec_\perp\cdot\textbf{grad}_{\vvec_\perp}f-F_z\fracp{f}{v_\zeta}=0\,,\label{vlasovscale}
\end{eqnarray}

\noindent Next, defining the quantity $\Ecal_\perp'=\p{E_x'-\beta B_y',E_y'+\beta B_x'}$, one easily verifies that $\Ecal_\perp=\bar{E}\Ecal_\perp'$. Accordingly, applying these dimensionless variables and dropping still the primes, Ampere's law (\ref{Amp11}-\ref{Amp12}) and the Poisson equation (\ref{GaussE1}) give
\bea
&&\eta\fracp{\Evec_\perp }{t }+\frac{1}{\beta}\fracp{}{\zeta }\p{\Ecal_\perp -\p{1-\beta^2}\Evec_\perp }-\curlvp B_z =-\eta\Jvec_\perp\,, \label{Ampexp1}\\
&&\eta\fracp{E_z }{t }+\frac{1}{\beta}\divp \p{\Ecal_\perp -\p{1-\beta^2}\Evec_\perp }=\eta J_\zeta \,,\label{Ampexp2}\\
&&\divp \Evec_\perp -\fracp{E_z }{\zeta }=\rho\,, \label{Ampexp3}
\eea
whereas Faraday's law (\ref{eq:eq2}-\ref{eq:eq3}) and the absence of monopoles equations (\ref{eq:eq4}) are written
\bea
&&\eta\fracp{\Bvec_\perp }{t }+\fracp{}{\zeta }\p{\Ecal_\perp \times\khat}+\curlvp E_z =\mathbf{0}\,, \label{Farexp1}\\
&&\eta\fracp{B_z }{t }+\curlsp \Ecal_\perp =0 \,,\label{Farexp2}\\
&&\divp \Bvec_\perp -\fracp{B_z }{\zeta }=0 \,.\label{Farexp3}
\eea
In the above equations, the right-hand sides $\rho$ and $(\Jvec_\perp,J_\zeta)$ fulfill the charge conservation equation
\be
 \label{chargeconsexp}
\eta\p{\fracp{\rho }{t }+\divp \Jvec_\perp +\fracp{J_\zeta }{\zeta }}=0\,.
\ee
Finally, the electromagnetic force $\Fvec=(\Fvec_\perp,F_z)$ takes the form
\bea
\Fvec_\perp &=&\Ecal_\perp +\eta\p{B_z \vvec_\perp +v_\zeta \Bvec_\perp }\times\khat \,,\label{Forceasymp1}\\
F_z &=&E_z +\eta\p{v_x B_y -v_y B_x }\,. \label{Forceasymp2}
\eea

\noindent We turn to the boundary conditions. The scaled electric field \textbf{E} obeys the same boundary conditions on the perfectly conducting boundary of the tube, together with the scaled analogous of (\ref{BC-Ecal}), i.e. $\Ecal_\perp\cdot\tauv=\beta \Bvec_\perp\cdot\nuv$. Concerning the scaled magnetic field $(\textbf{B}_\bot, B_z)$, we get from (\ref{eq:eq05}-\ref{eq:eq052})
\begin{eqnarray*}
 &&(\eta \frac{\partial}{\partial t}+\beta\frac{\partial}{\pa\zeta})\textbf{B}_\bot \cdot \nuv=0 ,
\qquad \eta\int_\Omega \partialderiv{B_z}{t}\diff\textbf{x}_{\bot} +\beta \oint_{\Gamma} \Bvec_\bot \cdot \nuv \diff l =0,
\qquad
\int_{\Omega}(\eta \frac{\partial}{\partial t}+\frac{1}{\eta}\frac{\partial}{\pa\zeta}) B_z \diff\textbf{x}_\bot=0,
\end{eqnarray*}
whereas, for $\textbf{x}_\bot \in \Omega ,\zeta=0$, we get $\textbf{E} =\textbf{0}\,, \textbf{B}=\textbf{B}^{e}$ and for $\textbf{x}_\bot \in \Omega ,\zeta=Z$, we obtain ${\Ecal}_{\bot}=0$.

\section{An asymptotic expansion}
\label{Asymptotic}

In order to derive a paraxial model, let us now rewrite the scaled Vlasov-Maxwell equations using expansions of
the quantities f, $\rho$, $\Jvec$, \textbf{E}, \textbf{B},  $\Ecal_\bot$ and \textbf{F}  in powers of the small parameter $\eta$, namely:\\

 $f=f^0+\eta f^1+\eta^2 f^2+...$, \qquad  $\rho=\rho^0+\eta \,\rho^1+\eta^2 \rho^2+...$, \qquad$\Jvec=\Jvec^0+\eta\, \Jvec^1+\eta^2 \Jvec^2+...$, \\

$\textbf{E}=\textbf{E}^0+\eta \textbf{E}^1+\eta^2 \textbf{E}^2+...$,
\, \quad $\!\textbf{B}=\textbf{B}^0+\eta \textbf{B}^1+\eta^2 \textbf{B}^2+...$,
\,  $\,\,\,\Ecal_\bot={\Ecal_\bot}^0+\eta {\Ecal_\bot}^1+\eta^2 {\Ecal_\bot}^2+...$,\\

\, $\textbf{F}=\textbf{F}^0+\eta \textbf{F}^1+\eta^2 \textbf{F}^2+...$ .\\

\noindent Then, we replace formally in the scaled Vlasov-Maxwell equations the functions by their asymptotic expansions, and we
identify the coefficients of $\eta^0$, $\eta^1$, etc. We begin by applying these expansions to the Vlasov equation (\ref{vlasovscale}). We get:
\begin{itemize}
\item at the zeroth order
$$
\fracp{f^0}{t}+\vvec_\perp\cdot\textbf{grad}_\perp f^0+v_\zeta\fracp{f^0}{\zeta}\
+\Fvec_\perp^0\cdot\gradv_{\vvec_\perp}f^{0}+
F_z^0 \,\fracp{f^{0}}{v_\zeta}=0\,,
$$
\item or at the first order
$$
\fracp{f^1}{t}+\vvec_\perp\cdot\textbf{grad}_\perp f^1+v_\zeta\fracp{f^1}{\zeta}+
\Fvec_\perp^0\cdot\gradv_{\vvec_\perp}f^{1}+\Fvec_\perp^1\cdot\gradv_{\vvec_\perp}f^{0}+
F_z^0\,\fracp{f^{1}}{v_\zeta}+F_z^1\,\fracp{f^{0}}{v_\zeta}=0\,.
$$
\end{itemize}
More generally, one can write out this equation for powers of $\eta$, that is, for $n^{{th}}$ order:
\be
\label{VlasAsymGen}
\fracp{f^n}{t}+\vvec_\perp\cdot\textbf{grad}_\perp f^n+v_\zeta\fracp{f^n}{\zeta}+\sum\limits_{i=0}^n\Fvec_\perp^i\cdot\gradv_{\vvec_\perp}f^{n-i}+\sum\limits_{i=0}^nF_z^i\,\,\fracp{f^{n-i}}{v_\zeta}=0
\ee
in which we use the convention that the negative superscripts vanish.\\

\noindent Hence, for determining the asymptotic expansion of the distribution function $f$ up to a given order $n$ in $\eta$, it is enough to know the expansion of the transverse and longitudinal electromagnetic force $\Fvec_\perp,F_z$  up to their $n$-th order. Then, using the expressions (\ref{Forceasymp1}-\ref{Forceasymp2}) of the forces, we get, with the same convention on the negative superscript:
\bea
\Fvec_\perp^n&=&\Ecal_\perp^n+\p{B_z^{n-1}\vvec_\perp+v_\zeta\Bvec_\perp^{n-1}}\times\khat\label{eq:eq27-1}\,,\\
F_z^n&=&E_z^n+\vvec_\perp\cdot\p{\Bvec_\perp^{n-1}\times\khat}\label{eq:eq27-2}\,.
\eea

\noindent In these conditions, the asymptotic expressions of these forces are entirely determined  as soon as we know  the expansions of $\Ecal_\perp$ and $E_z$ up to the $n$-th order and $\Evec_\perp$\footnote{$\Evec_\perp$ does not appear explicitly in the forces (\ref{eq:eq27-1}-\ref{eq:eq27-2}), but is required to compute $B_{z}$}, $\Bvec_\perp$ and $B_z$ up to the $\p{n-1}$-th order. 
Our aim now is to determine equations that characterize these ``required'' electromagnetic asymptotic fields.\\

\noindent For this purpose, we apply theses expansions to Maxwell's equations. Remark that all the terms where a time derivative is involved is multiplied by $\eta$, so they do not appear in the zeroth order. Hence, we obtain
\begin{itemize}
\item for Ampere's law and the Poisson equations (\ref{Ampexp1}-\ref{Ampexp3})
\begin{eqnarray*}
&&\fracp{}{\zeta}\p{\Ecal_\perp^0-\p{1-\beta^2}\Evec_\perp^0}-\beta\, \curlvp B_z^0=0\,,\\
&&\divp\p{\Ecal_\perp^0-\p{1-\beta^2}\Evec_\perp^0}=0\,,\\
&&\divp\Evec_\perp^0-\fracp{E_z^0}{\zeta}=\rho^0\,,
\end{eqnarray*}
whereas Faraday's law and the absence of monopole equations (\ref{Farexp1}-\ref{Farexp3}) yield
\begin{eqnarray*}
&&\fracp{}{\zeta}\p{\Ecal_\perp^0\times\khat}+\curlvp E_z^0=\mathbf{0}\,,\\
&&\curlsp\Ecal_\perp^0=0\,,\\
&&\divp\Bvec_\perp^0-\fracp{B_z^0}{\zeta}=0\,.
\end{eqnarray*}
Finally, the charge conservation equation (\ref{chargeconsexp}) leads to 
\begin{eqnarray*}
\fracp{\rho^0}{t}+\divp\Jvec_\perp^0+\fracp{J_\zeta^0}{\zeta}=0\,.
\end{eqnarray*}
\end{itemize}
On the contrary at the first order, the terms with a time derivative do appear, with an index $^0$. More precisely, we have, for Ampere's law
\begin{eqnarray*}
&&\fracp{\Evec_\perp^{0}}{t}+\frac{1}{\beta}\fracp{}{\zeta}\p{\Ecal_\perp^1-\p{1-\beta^2}\Evec_\perp^1}-\curlvp B_z^1=-\Jvec_\perp^{0}\,,\\
&&\fracp{E_z^{0}}{t}+\frac{1}{\beta}\divp\p{\Ecal_\perp^1-\p{1-\beta^2}\Evec_\perp^1}=J_\zeta^{0}\,,
\end{eqnarray*}
and for Faraday's law
\begin{eqnarray*}
&&\fracp{\Bvec_\perp^{0}}{t}+\fracp{}{\zeta}\p{\Ecal_\perp^1\times\khat}+\curlvp E_z^1=\mathbf{0}\,,\\
&&\fracp{B_z^{0}}{t}+\curlsp\Ecal_\perp^1=0\,.
\end{eqnarray*}
The other equations have the same expression simply by replacing index $0$ with index $1$. More generally, these expansions can be written out by the general following expressions for the n-th order. We obtain, for the electric field, still using the same convention on the negative superscript):
\bea
&&\fracp{\Evec_\perp^{n-1}}{t}+\frac{1}{\beta}\fracp{}{\zeta}\p{\Ecal_\perp^n-\p{1-\beta^2}\Evec_\perp^n}-\curlvp B_z^n=-\Jvec_\perp^{n-1}\label{eq:eq19}\,,\\
&&\fracp{E_z^{n-1}}{t}+\frac{1}{\beta}\divp\p{\Ecal_\perp^n-\p{1-\beta^2}\Evec_\perp^n}=J_\zeta^{n-1}\label{eq:eq9}\,,\\
&&\divp\Evec_\perp^n-\fracp{E_z^n}{\zeta}=\rho^n\label{eq:eq8}\,,
\eea
whereas, for the magnetic field, one gets
\bea
&&\fracp{\Bvec_\perp^{n-1}}{t}+\fracp{}{\zeta}\p{\Ecal_\perp^n\times\khat}+\curlvp E_z^n=0\,, \label{eq:eq11}\\
&&\fracp{B_z^{n-1}}{t}+\curlsp\Ecal_\perp^n=0\,, \label{eq:eq14}\\
&&\divp\Bvec_\perp^n-\fracp{B_z^n}{\zeta}=0\,, \label{eq:eq18}
\eea
and the charge conservation equation is expressed as
\bea
&&\fracp{\rho^n}{t}+\divp\Jvec_\perp^n+\fracp{J_\zeta^n}{\zeta}=0\,.
\eea

\noindent For the sake of completeness, we finally present the boundary conditions, that can be expressed, for $\xvec_\perp\in\Gamma,~\zeta\in\p{0,Z}$:
\begin{eqnarray}
&&\Evec_\perp^n\cdot\tauv=0,  \qquad \qquad E_z^n=0,  \qquad\qquad \Ecal_\perp^n\cdot\tauv=\beta\Bvec_\perp^n\cdot\nuv,\label{BCEadim}\\
&&\hspace*{-0.5cm}\p{\fracp{\Bvec_\perp^{n-1}}{t}+\beta\fracp{\Bvec_\perp^n}{\zeta}}\cdot\nuv=0, \,\int\limits_\Omega\fracp{B_z^{n-1}}{t}\diff\xvec_\perp+\beta\oint\limits_\Gamma\Bvec_\perp^n\cdot\nuv\diff l=0, \,
\int\limits_\Omega\p{\fracp{B_z^{n-1}}{t}+\beta\fracp{B_z^n}{\zeta}}\diff\xvec_\perp=0 \label{BCEadim}
\end{eqnarray}

 \noindent  As a consequence, one can write the following lemmas that characterize the different field component, at a given order $n$. One has first, for the longitudinal electric component $E_z^n$\\
 
\begin{Lem}
\label{LemEzn}
The $n$-th order component $E_z^n$ is the unique solution to
\be
\label{eq:eq13}
\left \{
\begin{array}{l}
\ds{\Delta_\perp}E_z^n+\p{1-\beta^2}\frac{\partial^2E_z^n}{\partial{\zeta}^2}=\\
\hspace*{3.cm} \fracp{}{t}\p{\beta\fracp{E_z^{n-1}}{\zeta}+\curlsp\Bvec_\perp^{n-1}}
\ds-\fracp{}{\zeta}\p{\beta J_\zeta^{n-1}+\p{1-\beta^2}\rho^n}
\textup{ in }\Omega\\
E_z^n=0  \textup{ on }  \Gamma
\end{array}
  \right.
\ee
\end{Lem}
\begin{prooff}
 Inserting (\ref{eq:eq8}) into (\ref{eq:eq9}) gives 
 \begin{equation}
 \label{InproofLezn}
 \divp\Ecal_\perp^n-\p{1-\beta^2}\fracp{E_z^n}{\zeta}=\p{1-\beta^2}\rho^n+\beta J_\zeta^{n-1}-\beta\fracp{E_z^{n-1}}{t}\,.
 \end{equation}
 Then, differentiating this relation with respect to $\zeta$ and adding the $\curlsp$ of (\ref{eq:eq11}) gives the desired result. 
\end{prooff}

\noindent Then, $E_z^n$ and quantities of the previous order $n-1$ are used to compute the pseudo-field $\Ecal_\perp^n$:

\begin{Lem}
\label{LemEcaln}
The $n$-th order component $\Ecal_\perp^n$ is the unique solution to
\begin{equation}
\label{eq:eq21}
\left \{
\begin{array}{l}
\ds\curlsp\Ecal_\perp^n=-\fracp{B_z^{n-1}}{t}\\
\ds\divp\Ecal_\perp^n=\p{1-\beta^2}\p{\fracp{E_z^n}{\zeta}+\rho^n}+\beta\p{J_\zeta^{n-1}-\fracp{E_z^{n-1}}{t}}
 \textup{in }\Omega\\
\ds\oint\limits_\Gamma\Ecal_\perp^n\cdot\tauv\diff l=-\int\limits_\Omega\fracp{B_z^{n-1}}{t}\diff\xvec_\perp 
\end{array}
 \right.
\end{equation}
\end{Lem}
\begin{prooff}
Since $E_z^n$ is known from (\ref{eq:eq13}), getting the equations is straightforward from (\ref{eq:eq14}) and (\ref{InproofLezn}). The boundary conditions are easily obtained from their expressions above. 
\end{prooff}

\noindent Similarly,  one gets the system that solves the transverse electric field $\Evec_\perp^n$, required after that to obtain the  transverse magnetic field $\Bvec_\perp^n$ (see below Lemma \ref{LemBperpn})
 \begin{Lem}
 \label{LemEperpn}
The $n$-th order component $\Evec_\perp^n$ is the solution to
\be
\left\{
\begin{array}{l}
\ds\curlvp\p{\curlsp\Evec_\perp^n}-\p{1-\beta^2}\frac{\partial^2\Evec_\perp^n}{\partial{\zeta}^2}\\
\hspace*{3.cm}\ds=-\frac{\partial^2\Ecal_\perp^n}{\partial{\zeta}^2}-\curlvp\p{\fracp{B_z^{n-1}}{t}}-\beta\fracp{}{\zeta}\p{\fracp{\Evec_\perp^{n-1}}{t}+\Jvec_\perp^{n-1}}
\quad\textup{ in }\,\Omega\,,\\
\divp\Evec_\perp^n=\fracp{E_z^n}{\zeta}+\rho^n \quad\textup{ in }\,\Omega\,,\\
\ds\Evec_\perp^n\cdot\tauv=0 \qquad \textup{ on } \Gamma\,,
\end{array}
\right.
\ee
\end{Lem}

\noindent \begin{prooff}
Computing $\curlsp\Ecal_\perp^n:=\curlsp(\Evec_\perp^n-\beta\Bvec_\perp^n\times\khat)$ and using (\ref{eq:eq14}-\ref{eq:eq18}) gives
\be
\fracp{B_z^n}{\zeta}=-\frac{1}{\beta}\fracp{B_z^{n-1}}{t}-\frac{1}{\beta}\curlsp\Evec_\perp^n
\ee
Combining it with the derivative of (\ref{eq:eq19}) with respect to $\zeta$ gives the result, $\Ecal_\perp^n$ being known from (\ref{eq:eq21}).
\end{prooff}

\noindent The two last results are concerned with the magnetic filed. First we have, for the transverse component:
\begin{Lem}
\label{LemBperpn}
\label{LemmaBerp}
The $n$-th order component $\Bvec_\perp^n$ is the unique solution to
\be
\label{eq:eq23}
\left\{
\begin{array}{l}
\ds\curlsp\Bvec_\perp^n=\fracp{E_z^{n-1}}{t}+\beta\divp\Evec_\perp^n-J_\zeta^{n-1}\\
\ds\divp\Bvec_\perp^n=-\frac{1}{\beta}\p{\curlsp\Evec_\perp^n+\fracp{B_z^{n-1}}{t}}
\textup{in }\Omega\\
\ds\oint\limits_\Gamma\Bvec_\perp^n\cdot\nuv\diff l=-\frac{1}{\beta}\int\limits_\Omega\fracp{B_z^{n-1}}{t}\diff\xvec_\perp \textup{on }\Gamma 
\end{array}
\right.
\ee
\end{Lem}
\begin{prooff}
Computing $\divp\Ecal_\perp^n:=\divp(\Evec_\perp^n-\beta\Bvec_\perp^n\times\khat)$ in combination with (\ref{eq:eq9}) gives one of the equations. The second one is obtained by combining (\ref{eq:eq22}) and (\ref{eq:eq14}), and the boundary condition is (\ref{BCEadim}). 
\end{prooff}

\noindent Finally, the longitudinal component $B_z^n$ is entirely determined by the magnetic field, and is characterized by:
\begin{Lem}
\label{LemBzn}
The $n$-th order component $B_z^n$ is the unique solution to
\be
\left\{
\begin{array}{l}
\ds\fracp{B_z^n}{\zeta}=\divp\Bvec_\perp^n \textup{in }\Omega\\
\ds \int\limits_\Omega\fracp{B_z^n}{\zeta}\diff\xvec_\perp=-\frac{1}{\beta}\int\limits_\Omega\fracp{B_z^{n-1}}{t}\diff\xvec_\perp 
\end{array}
\right.
\ee
\end{Lem}
\begin{prooff}
$\Bvec_\perp^n$ being known from (\ref{eq:eq23}), the equation is given by (\ref{eq:eq18}). The boundary condition is straightforward to obtain. 
\end{prooff} 
\section{The paraxial model}
\label{ParaxialModel}
\noindent We are now ready to introduce the paraxial model, which provides an approximation of the distribution function $f$ which is formally $n$ order accurate in $\eta$: this means that the asymptotic expansions of $f$ in the Vlasov-Maxwell and in the paraxial model coincide up the order $n$ in $\eta$. We will derive this model coming back to the physical variables, by using the scaling factors as introduced in Section \ref{Scaling}\footnote{remember that we dropped the $'$ in the previous section}.\\

\noindent To begin with, let us derive the equations satisfied by $E_z^n$. Assuming the knowledge of the data $(\rho, \Jvec)$, and of the fields up to the order $n-1$,  we obtain, from Lemma \ref{LemEzn}
\be
\label{eq:eq25}
\left\{
\begin{array}{l} 
\ds{\Delta_\perp}^2E_z^n+\p{1-\beta^2}\frac{\partial^2E_z^n}{\partial{\zeta}^2}=\frac{1}{c}\pp{\fracp{}{t}\p{\beta\fracp{E_z^{n-1}}{\zeta}+\curlsp c\Bvec_\perp^{n-1}}\right.\\
\hspace*{4.8cm}\ds\left.-\frac{1}{\varepsilon_0}\fracp{}{\zeta}\p{\beta J_\zeta^{n-1}+\p{1-\beta^2}c\rho^n}}
\,\,\, \textup{in } \Omega\\
E_z^n=0 \,\,\textup{on }\Gamma
\end{array}
\right.
\ee
Let us now deal with the transverse electric field.  From $E_z^n$, one can compute $\Ecal_\perp^n$ by solving to a quasi-static model, that is written, following Lemma \ref{LemEcaln}
\be
\left\{
\begin{array}{l}
\ds\curlsp\Ecal_\perp^n=-\fracp{B_z^{n-1}}{t}\\
\ds\divp\Ecal_\perp^n=\p{1-\beta^2}\p{\fracp{E_z^n}{\zeta}+\frac{\rho^n}{\varepsilon_0}}+\frac{\beta}{\varepsilon_0c}J_\zeta^{n-1}-\frac{\beta}{c}\fracp{E_z^{n-1}}{t}
\,\,\, \textup{in }\Omega\\
\ds\oint\limits_\Gamma\Ecal_\perp^n\cdot\tauv\diff l=-\int\limits_\Omega\fracp{B_z^{n-1}}{t}\diff\xvec_\perp 
\end{array}
\right.
\ee
In our paraxial model, even if $\Evec_\perp^n$ does not appear explicitly in the expression of the forces, there is yet a need to compute it as is required to obtain $B_{z}$. Following Lemma \ref{LemEperpn}, we have
\be
\left\{
\begin{array}{l}
\ds\curlvp\p{\curlsp\Evec_\perp^n}-\p{1-\beta^2}\frac{\partial^2\Evec_\perp^n}{\partial{\zeta}^2}\\
\hspace*{2.8cm}\ds=-\frac{\partial^2\Ecal_\perp^n}{\partial{\zeta}^2}-\curlvp\p{\fracp{B_z^{n-1}}{t}}-\frac{\beta}{c}\fracp{}{\zeta}\p{\fracp{\Evec_\perp^{n-1}}{t}+\frac{\Jvec_\perp^{n-1}}{\varepsilon_0}}
\,\,\, \textup{in }\Omega\\
\divp\Evec_\perp^n=\fracp{E_z^n}{\zeta}+\frac{\rho^n}{\varepsilon_0}\quad\textup{ in }\,\Omega\,,\\
\ds\Evec_\perp^n\cdot\tauv=0  \,\,\textup{on }\Gamma
\end{array}
\right.
\ee
This allows us to compute now the transverse magnetic field $\Bvec_\perp^n$, by solving, following Lemma \ref{LemBperpn}, the quasi-static system of equations
\be
\left\{
\begin{array}{l}
\ds\curlsp\Bvec_\perp^n=\frac{1}{c^2}\fracp{E_z^{n-1}}{t}+\frac{\beta}{c}\divp\Evec_\perp^n-\mu_0 J_\zeta^{n-1}\,\,\,\textup{in }\Omega\\\
\ds\divp\Bvec_\perp^n=-\frac{1}{\beta c}\p{\curlsp\Evec_\perp^n+\fracp{B_z^{n-1}}{t}}
\,\,\,\textup{in }\Omega\\
\ds\oint\limits_\Gamma\Bvec_\perp^n\cdot\nuv\diff l=-\frac{1}{\beta c}\int\limits_\Omega\fracp{B_z^{n-1}}{t}\diff\xvec_\perp \,\,\,\textup{on }\Gamma
\end{array}
\right.
\ee
Finally, one can obtain the longitudinal magnetic field of order $n$ by solving the simple equation, deduced from Lemma \ref{LemBzn}  
\be
\label{eq:eq30}
\left\{
\begin{array}{l}
\ds\fracp{B_z^n}{\zeta}=\divp\Bvec_\perp^n \,\,\,\textup{in }\Omega\\
\ds \int\limits_\Omega\fracp{B_z^n}{\zeta}\diff\xvec_\perp=-\frac{1}{\beta}\int\limits_\Omega\fracp{B_z^{n-1}}{t}\diff\xvec_\perp 
\end{array}
\right.
\ee

\noindent The paraxial model proposed here is hierarchical and closed for each order: the zeroth order fields allow to solve the first order etc. Note also that the time derivatives being on the left-hand side, the model is quasi-static and not time-dependent. In addition, the $n$-th order fields are required only for $\Ecal_\perp$ and $E_z$, whereas it is sufficient to know the other fields up to the $\p{n-1}$-th order.\\

\noindent We can summarize our main result in the following theorem:\\

 \begin{theorem}
Equations (\ref{eq:eq25}-\ref{eq:eq30}) determine the triple $(\Evec^{n}, \Bvec^{n},\Ecal_\perp^n)$ from the data $(\rho, \Jvec)$, and $(\Evec^{l}, \Bvec^{l},\Ecal_\perp^{l})$, for $0 \leq l \leq n-1$,  
in a unique way. Moreover, the paraxial model provides an approximation of the distribution function $f$ which is formally of order $n$ accurate in $\eta$, namely, the asymptotic expansions of $f$ in the Vlasov-Maxwell and in the paraxial model coincide up the $n$ order in $\eta$.\\
\end{theorem}

\section{Conclusion}

\noindent In this Note, we proposed a new family of paraxial asymptotic models which approximate the non-relativistic Vlasov-Maxwell equations. It has been derived by introducing a small parameter $\eta=\ds\frac{\overline{v}}{c}$, and is $n$-th order accurate, for $n \in \N$. In these conditions, one can easily choose the complexity of the model one wants to use, depending on the required accuracy. In addition, this family of models is simpler than the Vlasov-Maxwell equations - for instance they are not time-dependent but only static or quasi-static - which allows us to implement simple and efficient numerical schemes, like particle-in-cell techniques. Hence, this approach would be very powerful in its ability to get fast and easy to implement algorithms.


%

\begin{thebibliography}{00}
%
\bibitem{weld} W.J. Harris, U.S. Patent No. 3,271,556. Washington, DC: U.S. Patent and Trademark Office, 1966.
%
\bibitem{lith} M.J. Madou, Manufacturing techniques for microfabrication and nanotechnology, 2, CRC press, 2011.
%
\bibitem{fusion} R.B. Miller, An introduction to the physics of intense charged particle beams, Springer, 1984.
%
\bibitem{gyrotron} B. Danly, G. Bekefi, R. Davidson, R. Temkin, T. Tran, J. Wurtele, Principles of gyrotron powered electromagnetic wigglers for free-electron lasers. IEEE journal of quantum electronics, 23(1), 103--116, 1987.
%
\bibitem{tran} T.M. Tran, J.S. Wurtele, Free-electron laser simulation techniques, Physics Reports, 195(1), 1--21, 1990.
%
\bibitem{Laws88} J.D. Lawson, The Physics of Charged Particle Beams, Oxford, Clarendon Press, 1988.
%
\bibitem{reiser2008book} M. Reiser, Theory and design of charged particle beams, John Wiley \& Sons, 2008.
%
\bibitem{vlasov} A. Vlasov, On the kinetic theory of an assembly of particles with collective interaction, Russ. Phys. J., 9, 25--40, 1945.
%
\bibitem {BiLa85} C.K. Birdsall and A.B. Langdon, Plasmas Physics via Computer Simulation, (New York: Mac.Graw-Hill, 1985).
%
\bibitem {AsCL18} F. Assous, P. Ciarlet, Jr., S. Labrunie, Mathematical Foundations of Computational Electromagnetism, Appl. Math. Sc., AMS 198, Springer, 2018.
%
\bibitem {ADHRS93} F. Assous, P. Degond, E. Heintz\'e, P.A. Raviart, J. Segr\'e, On a finite element method for solving the three dimensional Maxwell equations,  {J. Comput. Physics}, 109(2), 222--237, 1993.
%
\bibitem{AsDS92} F. Assous, P. Degond, J. Segr\'e, A particle-tracking method for 3D electromagnetic PIC codes on unstructured meshes,  {Comput. Phys. Comm.} { 72} 105--114, 1992.
%
\bibitem{DeRa93} P. Degond,P.-A. Raviart, On the paraxial approximation of the stationary Vlasov-Maxwell, {Math. Mod. Meth. Appl. Sci.}, { 3(4)} 513--562, 1993.
%
 \bibitem{LaMR94} G. Laval, S. Mas-Gallic, P.-A. Raviart, Paraxial approximation of ultrarelativistic intense beams, {Numer. Math.}, {69(1)}, 33--60, 1994.
%
\bibitem{MoMW96} M.A. Mostrom, D. Mitrovich, D.R. Welch,The ARCTIC Charged Particle Beam Propagation Code,
{J. Comput. Physics}, {128(2)} 489--497, 1996.
%
\bibitem{RaSo96} P.A. Raviart, E. Sonnendr\"ucker, A hierarchy of approximate models for the Maxwell equations, {Numer. Math.}, {73(3)}, 329--372, 1996.
%
\bibitem{boyd} J.K. Boyd, E.P. Lee, S.S. Yu,  Aspects of three field approximations: Darwin, frozen, EMPULSE (No. UCID-20453). Lawrence Livermore National Lab., CA, USA, 1985.
%
\bibitem{slinker} S. Slinker, G. Joyce, J. Krall, R.F. Hubbard, ELBA a three dimensional particle simulation code for high current beams. In Proc. of the 14th Inter. Conf. Numer. Simul. Plasmas, Annapolis, 1991.
%
\bibitem{nouri} A. Nouri, Paraxial approximation of the Vlasov-Maxwell system: laminar beams, Mathematical Models and Methods in Applied Sciences, 4(02), 203--221,1994.
%
\bibitem{AsCh12} F. Assous, J. Chaskalovic, A New Paraxial Asymptotic Model for the Relativistic Vlasov-Maxwell Equations,
{C. R. Mecan. Acad. Sci}, {340}, 706-714, 2012.
%
\bibitem{AsTs09} F. Assous, F. Tsipis, Numerical paraxial approximation for highly relativistic beams, {Comput. Phys. Comm.}, {180-7}, 10861--1097, 2009.
%
\bibitem{AsCh11} F. Assous, J. Chaskalovic, Data mining techniques for scientific computing:  Application to asymptotic paraxial
approximations to model ultra relativistic particles, {J. Comput. Phys.}, {230}, 4811--4827, 2011.
%
\bibitem{AsCh14} F. Assous, J. Chaskalovic, A Paraxial Asymptotic Model for the Coupled Vlasov-Maxwell Problem in Electromagnetics, {\em J. Comput. Appl. Maths}, 270, 369--385, 2014.
%
\end{thebibliography}
\end{document}